\documentclass[11pt]{article}
\usepackage{amsmath}
\usepackage{amssymb}
\usepackage{amsthm}
\usepackage{amsfonts}
\usepackage[T1]{fontenc}
\usepackage[latin1]{inputenc}

\usepackage{epsfig}
\usepackage[lflt]{floatflt}
\usepackage{graphicx}
\usepackage{subfigure}

\usepackage{array}
\usepackage{tabularx}
\usepackage{graphicx}
\usepackage{verbatim}

\usepackage{fancyhdr}
\usepackage{color}
\usepackage[affil-sl]{authblk}
\usepackage{cancel}
\usepackage[english]{babel}
\usepackage{fancyhdr}
\usepackage{makeidx}
\usepackage{latexsym}
\usepackage{geometry}

\newtheorem{note}{Note}
\newtheorem{theorem}{Theorem}
\newtheorem{definition}{Definition}
\newtheorem{proposition}{Proposition}
\newtheorem{lemma}{Lemma}
\newtheorem{corollary}{Corollary}
\newcommand{\CC}{{\mathcal C}}
\newcommand{\WW}{{\mathcal W}}
\newcommand{\MM}{{\mathcal M}}
\newcommand{\KK}{{\mathcal K}}

\newcommand{\bbR}{{\mathbb R}}
\newcommand{\bbN}{{\mathbb N}}

\newcommand{\e}{\epsilon}
\newcommand{\al}{\alpha}
\newcommand{\la}{\lambda}
\newcommand{\p}{\partial}
\newcommand{\be}{\beta}
\newcommand{\G}{\Gamma}
\def\blacksquare{\hbox{\vrule width 4pt height 4pt depth 0pt}}
\def\qed{\ \ \ \hbox{}\nolinebreak\hfill $\blacksquare \  \  \  \  $ \par{}\medskip}
\begin{document}
%\begin{center}
\LARGE
\textbf{On the Initial--Boundary--Value Problem for the Time--Fractional Diffusion Equation in the Quarter Plane}
%\end{center}
\medskip

\normalsize
\begin{center}
D. Goos$^{1}$, G. Reyero$^{2}$,  S. Roscani$^3$ and E. Santillan Marcus$^4$\\
\medskip
\small
$^{1,2,3,4}$Departamento de Matem\'{a}tica,  FCEIA, Universidad Nacional de Rosario, Pellegrini 250, Rosario, Argentina \\
\textcolor{blue}{demian@fceia.unr.edu.ar, greyero@fceia.unr.edu.ar, sabrina@fceia.unr.edu.ar, edus@fceia.unr.edu.ar} \\

$^4$ Departamento de Matem\'atica,
FCE, Universidad Austral, Paraguay 1950, Rosario, Argentina,\\
\textcolor{blue}{esantillan@austral.edu.ar}\\

$^{1,3}$ CONICET, Argentina.

\end{center}

%\thanks{}

\begin{abstract}
Taking into account the asymptotic behavior of some Wright functions and the existence of bounds for the Mainardi and the Wright function $W(-x,\frac{\alpha}{2}, 1)$ in $\mathbb{R}^+$ , three different  initial--boundary--value problems for the time--fractional diffusion equation in the quarter plane, where the time--fractional derivative is taken in the Caputo sense of order $\al$ $\in (0,1)$ are solved. Moreover, the limit   when $\al \nearrow 1$ of the respective solutions are analyzed, recovering  the respective solutions of the classical boundary--value problems when $\al=1$ and the fractional diffusion equation becomes the heat equation.
\end{abstract}
\small
\noindent \textbf{2000 Mathematics Subject Classiffcation: }Primary: 26A33, 35S10, 35S15. Secondary: 33E20, 35C15. \\
\textbf{Keywords: } Time--fractional diffusion equation; Caputo derivative;Mainardi function; integral representation of solution.

\normalsize

\section{Introduction}

\noindent The one--dimensional Heat Equation has become the paradigm for the all--embracing study of parabolic partial differential equations, linear and nonlinear. A methodical
development of a variety of aspects of this paradigm can be seen in \cite{Cannon, Friedman, Widder}. \\
\noindent 
This paper deals with two problems associated to the time--fractional diffusion equation, obtained from the standard diffusion equation by replacing the first order time--derivative by a fractional derivative of order $\alpha > 0 $ in the Caputo sense:
\begin{equation}\label{FDE}   _0 D^{\alpha}_t u(x,t)=\lambda^2\,  u_{xx}(x,t), \quad  -\infty<x<\infty, \ t>0, \ 0<\al<1, \end{equation}
 where the fractional derivative in the Caputo sense of arbitrary order $\al >0$ is given by
$$\,_{a} D^{\alpha}f(t)=\begin{cases} \frac{1}{\Gamma(n-\al)}\int^{t}_{a}(t-\tau)^{n-\al-1} f^{(n)}(\tau)d\tau, &  n-1<\al<n\\
f^{(n)}(t), &   \al=n. \end{cases}$$
where  $n \in \bbN$ and $\G$ is the Gamma function defined by $\G(x)=
\int_0^\infty  w^{x-1}e^{-w}dw$.\\

The interest on equation (\ref{FDE}) has been in constant increase during the last 30 years. So many authors have studied it  \cite{ Diethelm, Fujita:1989,  Kilbas, Luchko,  FM-TheFundamentalSolution,  Podlubny, Wyss:1986} and, among the several applications that have been studied, Mainardi \cite{FM-libro} focused on the application to the theory of linear viscoelasticity. \\
A comprehensive analysis of the Cauchy problem associated to this equation can be found in \cite{EidelmanKochubei:2004} and a physical meaning  is  discussed in \cite{Garra} .\\

The two initial--boundary--value problems considered are:

\begin{equation}{\label{PVI-IC}}
\left\{\begin{array}{lll}
         _0 D^{\al}_t c(x,t)=\lambda^2\dfrac{\partial^2 c }{\partial x^2}(x,t) &  & 0<x<\infty, \, 0<t<T, \, 0<\al<1, \\

          c(x,0)=f(x) & \quad  \quad  & 0<x<\infty, \\

          c(0,t)=g(t) & \quad  \quad  & 0<t<T, \\

\end{array}\right.\end{equation}

\noindent and
\begin{equation}{\label{PVIF-IC}}
\left\{\begin{array}{lll}
         _0 D^{\al}_t c(x,t)=\lambda^2\dfrac{\partial^2 c}{\partial x^2}(x,t) &  & 0<x<\infty, \, 0<t<T, \, 0<\al<1, \\

          c(x,0)=f(x) & \quad  \quad  & 0<x<\infty, \\

          \dfrac{\partial c}{\partial x}(0,t)=g(t) & \quad  \quad  & 0<t<T. \\

           %\displaystyle\lim_{x\rightarrow +\infty}c(x,t)=0 &\quad  \quad & t>0
                                       \end{array}\right.\end{equation}

\noindent Some variants of this problems have been solved in \cite{ GLM:2000} and \cite{Wyss:1986}. The former using scale--invariant techniques and Laplace transform, and the last one in terms of Fox functions using Mellin transform. We propose here a different approach involving convolutions and we prove in each case that the function proposed is a solution of the considered problem.\\
   
\noindent The paper is presented as follows: Some useful properties  about the behavior of Wright functions are given in Section 1.
 In sections 2, 3 and 4 the two problems enunciated previously will be solved. At the end of sections 2 and 4 the limit   when $\al \nearrow 1$ of the respective solutions will be done, recovering  the respective solutions of the classical boundary--value problems when $\al=1$ and equation (\ref{FDE}) becomes the heat equation.\\

\section{Preliminaries. Some Results about the special functions involved.}
\begin{definition}\label{def W} For every $z\in \mathbb{C}$ , $\al>-1$ and $\be\in \bbR$ the Wright function is defined by
\begin{equation}\label{defi W} \WW(z;\al;\be)=\sum^{\infty}_{k=0}\frac{z^{k}}{k!\G(\al k+\be)} .\end{equation}
\end{definition}

\begin{definition}\label{def M} For every $z\in \mathbb{C}$ , $0<\nu<1$  the Mainardi function is defined by
\begin{equation}\label{defi M} \MM_{\nu}(z)=\WW(z;-\nu;1-\nu)=\sum^{\infty}_{k=0}\frac{z^{k}}{k!\G(-\nu k +1-\nu)} .\end{equation}
\end{definition}

\begin{note}\label{deriv W} This  series are absolutely convergent over compact sets and so its derivatives are easy to calculate:
$$ \frac{d}{dz}\WW(z,\al,\be)=\sum^{\infty}_{n=0}\frac{d}{dz}\frac{(z)^n}{n! \G\left(\al n+\be \right)}=\sum^{\infty}_{n=1}\frac{(z)^{n-1}}{(n-1)! \G\left(\al n+\be \right)}= \sum^{\infty}_{n=0}\frac{(z)^{n}}{n! \G\left(\al (n+1)+\be \right)}=$$
$$
\quad \quad \quad \quad\quad \quad = \sum^{\infty}_{n=0}\frac{(z)^{n}}{n! \G\left(\al n+(\al+\be) \right)}=\WW(z,\al,\al +\be).
\quad \quad\quad \quad \, \,  \quad \quad \quad\quad \quad
\quad \quad\quad \quad $$
For the special case of the Mainardi function, we have 
$$ \frac{d}{dz} \MM_{\nu}(z)=\WW(z;-\nu;1-2\nu) $$
\end{note}

\subsection{\textbf{Asymptotic  behavior.}}

\noindent  The following asymptotic  behavior for the Mainardi function was proved in \cite{On a Special function}.

     \begin{equation} {\label{M-asintotica}}\MM_\nu\left(\frac{x}{\nu}\right)\sim a(\nu) x^{(\nu - 1/2)/(1-\nu)}exp\left[ -b(\nu)x^{1/(1-\nu)} \right],
     \end{equation}

     $$\text{ where }\quad  a(\nu)=\frac{1}{\sqrt{2\pi(1-\nu)}}>0 \text{ and } \quad \quad b(\nu)=\frac{1-\nu}{\nu}\, . $$

\begin{theorem}\label{Teo de Wright} If $-1<\rho<0$, $y=-z$, $|arg\, y|\leq \min\{\frac{3}{2}\pi(1+\rho), \pi \}-\e $, $\e >0$, then
$$\WW(z,\rho, \be)=I(Y),  \quad Y\rightarrow \infty , $$
where $$I(Y)=Y^{1/2-\be}e^{-Y}\left\{ \sum_{m=0}^{M-1}A_mY^{-m}+O(Y^{-M}) \right\} \quad \text{ and }\quad  Y=(1+\rho)\left((-\rho)^{-\rho}y\right)^{\frac{1}{1+\rho}}. $$
The coefficients $A_m, \, m=0,1,...$ are defined by the asymptotic expansion
$$\hspace{-4cm} \frac{\G(1-\be-\rho t)}{2\pi(-\rho)^{-\rho t}(1+\rho)^{(1+\rho)(t+1)}\G(t+1)}=\sum\limits_{m=0}^{M-1}\frac{(-1)^m A_m}{\G\left((1+\rho)t+\be+\frac{1}{2}+m\right)}+$$
$$\hspace{9.5cm} +O\left( \frac{1}{\G\left((1+\rho)t+\be+\frac{1}{2}+M\right)}\right),  $$

 \noindent valid for $\arg t$, $\arg (-\rho t)$ and $\arg (1-\be-\rho t)$ all lying between $-\pi $ and $\pi $ and $t$ tending to infinity.

\end{theorem}

\noindent This theorem was proved in \cite{Wright}. The next results follows.

\begin{corollary}\label{lim W(-x,-al/2,1)=0}
$$\lim_{x \rightarrow \infty}\WW\left(-x,-\frac{\al}{2},1\right)=0
\quad \text{ and } \quad \lim_{x \rightarrow \infty}\WW\left(-x,-\frac{\al}{2},1+\frac{\al}{2}\right)=0. \quad  $$
\end{corollary}
\begin{corollary}\label{lim M(x)=0}
$$\lim_{x \rightarrow \infty}\MM_{\al/2}\left(x\right)=0. $$
\end{corollary}

\begin{corollary}\label{acotacion en el inf W(-s,-al/2,1-al)}
 If  $\, 0<\al<1$ and  $x \in \bbR^+$,  there exists $R>0$ such that,
 $$\quad \left|\WW(-x,-\frac{\al}{2},1-\al)\right|< P\left(bx^{\frac{1}{1-\frac{\al}{2}}}\right)\exp\left\{-bx^{\frac{1}{1-\frac{\al}{2}}}\right\}\, \quad  \forall \, x>R,$$
 where $P(x)$ is a polynomial function of degree less or equal than $ 1$  and $b=(1-\frac{\al}{2})(\frac{\al}{2})^{\frac{\al}{2-\al}}>0$.

\end{corollary}

\noindent \textit{Proof.}\\
Let us consider the function $\WW(-x,-\frac{\al}{2},1-\al)$.
$$ z=-x\Rightarrow y=x \quad \text{ and }\quad  arg\, y=0 \, .$$
By Theorem \ref{Teo de Wright}, taking
$$ Y=\left(1-\frac{\al}{2}\right)\left(\left(\frac{\al}{2}\right)^{\frac{\al}{2}}x\right)^{\frac{1}{1-\frac{\al}{2}}}=bx^{\frac{1}{1-\frac{\al}{2}}}\quad \quad \text{ and } \quad b=\left(1-\frac{\al}{2}\right)\left(\frac{\al}{2}\right)^{\frac{\al}{2-\al}}>0, $$

$$\WW\left(-x,-\frac{\al}{2},1-\al\right)=(bx^{\frac{1}{1-\frac{\al}{2}}})^{\al-1/2}\exp\{-bx^{\frac{1}{1-\frac{\al}{2}}}\}\left\{ \sum_{m=0}^{M-1}A_m(bx^{\frac{1}{1-\frac{\al}{2}}})^{-m}+O\left((bx^{\frac{1}{1-\frac{\al}{2}}})^{-M}\right) \right\}.  $$
Or equivalently,

$$ \frac{\WW(-x,-\frac{\al}{2},1-\al)}{(bx^{\frac{1}{1-\frac{\al}{2}}})^{\al-1/2}\exp\{-bx^{\frac{1}{1-\frac{\al}{2}}}\}}-
\sum_{m=0}^{M-1}A_m(bx^{\frac{1}{1-\frac{\al}{2}}})^{-m}=O\left((bx^{\frac{1}{1-\frac{\al}{2}}})^{-M}\right).
$$

%\noindent This means that there exists $R>0$ such that
%
%$$ \left|\frac{\frac{\WW(-x,-\frac{\al}{2},1-\al)}{(bx^{\frac{1}{1-\frac{\al}{2}}})^{\al-1/2}\exp\{-bx^{\frac{1}{1-\frac{\al}{2}}}\}}-
%\sum\limits_{m=0}^{M-1}A_m(bx^{\frac{1}{1-\frac{\al}{2}}})^{-m}}{(bx^{\frac{1}{1-\frac{\al}{2}}})^{-M}} \right|\leq C, \quad \text{ if  } x>R.  $$

%\noindent This expression is valid for all  $M\geq 1 $. In particular for  $M=1$,

\noindent Taking $M=1$, there exists $R>0$ such that
$$ \left|\frac{\frac{\WW(-x,-\frac{\al}{2},1-\al)}{(bx^{\frac{1}{1-\frac{\al}{2}}})^{\al-1/2}
\exp\{-bx^{\frac{1}{1-\frac{\al}{2}}}\}}-
A_0}{(bx^{\frac{1}{1-\frac{\al}{2}}})^{-1}} \right|\leq C \quad \text{ if  } x>R. $$
%\noindent that is to say\\
%$ (bx^{\frac{1}{1-\frac{\al}{2}}})^{\al-1/2}\left(  -C (bx^{\frac{1}{1-\frac{\al}{2}}})^{-1} +A_0 \right)\leq
%\frac{\WW(-x,-\frac{\al}{2},1-\al)}{
%\exp\{-bx^{\frac{1}{1-\frac{\al}{2}}}\}} \leq
%(bx^{\frac{1}{1-\frac{\al}{2}}})^{\al-1/2}\left(  C (bx^{\frac{1}{1-\frac{\al}{2}}})^{-1} +A_0 \right).$\\
\noindent Then
  \begin{equation}\label{cota 1 para W(-x,-al-2,1-al)}\left| \WW\left(-x,-\frac{\al}{2},1-\al\right) \right|\leq (bx^{\frac{1}{1-\frac{\al}{2}}})^{\al-1/2}\left(  C (bx^{\frac{1}{1-\frac{\al}{2}}})^{-1} + |A_0| \right)\exp\{-bx^{\frac{1}{1-\frac{\al}{2}}}\},  \quad \text{ if } x>R. \end{equation}

\noindent If $0<\al<\frac{1}{2}$, $(bx^{\frac{1}{1-\frac{\al}{2}}})^{\al-1/2}<\frac{1}{(bR^{\frac{1}{1-\frac{\al}{2}}})^{1/2-\al}}$.\\
If $\al=\frac{1}{2}$, $(bx^{\frac{1}{1-\frac{\al}{2}}})^{\al-1/2}=1$.\\
If    $\frac{1}{2}<\al<1$,  taking R large enough so that $bx^{\frac{1}{1-\frac{\al}{2}}}>1$ if $x>R$, it follows that
    $$ (bx^{\frac{1}{1-\frac{\al}{2}}})^{\al-1/2}<
    (bx^{\frac{1}{1-\frac{\al}{2}}}) .$$

\noindent Hence there exists two constants $B_0$ and $B_1$ depending on  $\al$ such that
$$ (bx^{\frac{1}{1-\frac{\al}{2}}})^{\al-1/2}<B_0+B_1
    (bx^{\frac{1}{1-\frac{\al}{2}}}). $$

\noindent Finally,

$$ (bx^{\frac{1}{1-\frac{\al}{2}}})^{\al-1/2}
\left(  C (bx^{\frac{1}{1-\frac{\al}{2}}})^{-1} +|A_0| \right)<\left(B_0+B_1(bx^{\frac{1}{1-\frac{\al}{2}}})\right)\left(  C (bx^{\frac{1}{1-\frac{\al}{2}}})^{-1} +|A_0| \right)=$$

$$=P_1\left((bx^{\frac{1}{1-\frac{\al}{2}}})^{-1}\right)+P_2\left(bx^{\frac{1}{1-\frac{\al}{2}}}\right)
\leq \tilde{C} +P_2\left(bx^{\frac{1}{1-\frac{\al}{2}}}\right)=P\left(bx^{\frac{1}{1-\frac{\al}{2}}}\right),  $$
where $P$ is a polynomial function of degree less or equal than $1$.

\noindent Therefore
$$  \quad \left|\WW(-x,-\frac{\al}{2},1-\al)\right|< P\left(bx^{\frac{1}{1-\frac{\al}{2}}}\right)\exp\{-bx^{\frac{1}{1-\frac{\al}{2}}}\}\, ,  \quad \text{ if } x>R.$$
\qed

\begin{corollary}\label{acotacion en el inf W(-x,-al/2,1)} If $0<\al<1$ and $x\in \bbR^+$ , there exists $R>0$ such that,
 $$\quad \left|\WW\left(-x,-\frac{\al}{2},1\right)\right|<Ke^{-bx}\, ,   \quad  \forall \, x>R\quad \,\text{ where } b=\left(1-\frac{\al}{2}\right)\left(\frac{\al}{2}\right)^{\frac{\al}{2-\al}}.$$

\end{corollary}

\bigskip

\subsection{\textbf{Some bounds and convergence. \\}}

\noindent The assertions in this subsection were proved in \cite{RoSa1}.

\begin{lemma}\label{M pos y decrec} If  $\,  0<\al<1$, $\MM_{\al/2}(x) $ is a strictly decreasing positive function in $\bbR^+$ .
\end{lemma}

\begin{corollary}\label{acotacion M} If $x>0$, $\MM_{\al/2}(x)<\frac{1}{\G\left(1-\frac{\al}{2}\right)}$  .
\end{corollary}

\begin{corollary}\label{acotacion W(-x,-al/2,1)} If  $\,  0<\al<1$ , $\WW\left(-x,-\frac{\al}{2},1\right)$ is a positive and decreasing function in  $\bbR^+$ such that \\$0<\WW\left(-x,-\frac{\al}{2},1\right)\leq 1, \, \forall \, x \in \bbR^+_0$.

\end{corollary}

\begin{note} Note that

$$ \hspace{-2cm}\WW\left(-x,-\frac{1}{2},1\right)= \int_{\infty}^x \left( \frac{\p}{\p x} \WW\left(-\xi,-\frac{1}{2},1\right)\right)d\xi=\int_{\infty}^x -\WW\left(-\xi,-\frac{1}{2},\frac{1}{2}\right)d\xi =$$
 $$\hspace{-1.4cm}=  \int_x^{\infty} \WW\left(-\xi,-\frac{1}{2},\frac{1}{2}\right)d\xi=\int_x^{\infty} \frac{1}{\sqrt{\pi}}e^{-\xi^2/4}d\xi=$$
$$\hspace{-3.4cm} =  \frac{2}{\sqrt{\pi}}\int_{x/2}^{\infty} \frac{1}{\sqrt{\pi}}e^{-\xi^2}d\xi = \,erfc\left(\frac{x}{2}\right). $$
\noindent Hence
$$ 1-\WW\left(-x,-\frac{1}{2},1\right)=erf\left(\frac{x}{2}\right) . $$

\end{note}

\begin{lemma}\label{conv M al/2 cuando al tiende a 1}If $x\in \bbR^+_0$ and $\al \in (0,1)$,
\begin{enumerate}
\item $\lim\limits_{\al\nearrow 1}\MM_{\al/2}\left(x\right)=\MM_{1/2}(x)=\frac{e^{-\frac{x^2}{4}}}{\sqrt{\pi}}$
\item $\lim\limits_{\al\nearrow 1}\left[1-\WW\left(-x,-\frac{\al}{2},1\right)\right]=\frac{1}{\sqrt{\pi}} erf\left(\frac{x}{2}\right) .$

\end{enumerate}

\end{lemma}

\bigskip

\section{Solving the Initial--Boundary--Value Problem for the Time--Fractional Diffusion Equation in the Quarter Plane with  Temperature--Boundary Condition.}

\bigskip

\noindent Let us consider problem (\ref{PVI-IC}). The principle of superposition  is valid due to the linearity of the Caputo derivative. Then, solve problem (\ref{PVI-IC}) is equivalent to solve the two auxiliary problems:
\begin{equation}{\label{PVI-ICA}}
\left\{\begin{array}{lll}
          _0D^{\al}_t c_1(x,t)=\lambda^2\dfrac{\partial^2 c_1}{\partial x^2}(x,t) &  & 0<x<\infty, \, 0<t<T, \, 0<\al<1, \\
          c_1(x,0)=f(x) & \quad  \quad  & 0<x<\infty. \\
          c_1(0,t)=0 & \quad  \quad  & 0<t<T, \\
 \end{array}\right.\end{equation}

\noindent and

\begin{equation}{\label{PVI-ICB}}
\left\{\begin{array}{lll}
          _0D^{\al}_t c_2(x,t)=\lambda^2\dfrac{\partial^2 c_2}{\partial x^2} (x,t) &  & 0<x<\infty, \, 0<t<T, \, 0<\al<1, \\
          c_2(x,0)=0 & \quad  \quad  & 0<x<\infty, \\
          c_2(0,t)=g(t) & \quad  \quad  & 0<t<T. \\
                                       \end{array}\right.\end{equation}

\noindent  Problem (\ref{PVI-ICA}) was solved  in  \cite{Liu-Xu} and its solution is given by

 \begin{equation}\label{c_1(x,t)}
 c_1(x,t)=\int^{\infty}_{0}\frac{1}{2 \lambda t^{\al/2}}\left[ \MM_{\al/2}\left(\frac{|x-\xi|}{t^{\al/2}}\right) - \MM_{\al/2}\left(\frac{|x+\xi|}{t^{\al/2}}\right) \right]f(\xi) d\xi\, 
 \end{equation}

\noindent where the function $\MM_{\al/2}(\cdot)$ is the Mainardi function defined in (\ref{defi M}) and $f$ is a continuous bounded function in $\bbR_0^+$ ( which guarantees that $c_1$ is a solution, see the Cauchy problem in \cite{Eidelman-Kochubei-LIBRO}).

\noindent  In \cite{RoSa1} it was proved that 

\begin{equation}{\label{z(x,t)}}
z(x,t)=A+ (B-A)\WW\left(\frac{-x}{\lambda t^{\al/2}},-\frac{\al}{2},1\right), 
 \end{equation}
 where $\WW\left(\cdot,-\frac{\al}{2},1\right)$ is the Wright function of parameters $-\frac{\al}{2}$ and $1$, defined in (\ref{defi W}), is a solution to problem

\begin{equation}{\label{PVI-A,B}}
\left\{\begin{array}{lll}
          _0D^{\al}_t z(x,t)=\lambda^2\dfrac{\partial^2z }{\partial x^2}(x,t) &  & 0<x<\infty, \, 0<t<T, \, 0<\al<1, \\

          z(x,0)=A & \quad  \quad  & 0<x<\infty, \\
 
          z(0,t)=B & \quad  \quad  & 0<t<T. \\
\end{array}\right.\end{equation}

\noindent Then, we can assure that 
\begin{equation}{\label{v(x,t)}}
v(x,t)=\WW\left(\frac{-x}{\lambda t^{\al/2}},-\frac{\al}{2},1\right)
 \end{equation}
 is a solution to problem
\begin{equation}{\label{PVI1}}
\left\{\begin{array}{lll}
          _0D^{\al}_t v(x,t)=\lambda^2\dfrac{\partial^2v }{\partial x^2}(x,t) &  & 0<x<\infty, \, 0<t<T, \, 0<\al<1, \\
          v(x,0)=0 & \quad  \quad  & 0<x<\infty, \\

          v(0,t)=1 & \quad  \quad  & 0<t<T. 
\end{array}\right.\end{equation}

\noindent  Taking into account  Note \ref{deriv W} and Corollary \ref{lim M(x)=0}, function $v$ can be expressed as
  $$ v(x,t)=\int_0^t\frac{\p}{\p \tau}\WW\left(\frac{-x}{\lambda \tau^{\al/2}},-\frac{\al}{2},1\right)d\tau= \int_0^t \MM_{\al/2}\left(\frac{x}{\lambda \tau^{\al/2}}\right)\frac{\al x}{2\lambda \tau^{\al/2+1}}d\tau. $$	

\noindent Let be 
%$$  K(x,\tau)=\begin{cases} M_{\al/2}\left(\frac{x}{\lambda \tau^{\al/2}}\right)\frac{x}{\lambda \tau^{\al/2+1}}\frac{\al}{2} & 0<\tau<t\\
%0 & \text{ else }. \end{cases}$$

$$  K(x,t)= \MM_{\al/2}\left(\frac{x}{\lambda t^{\al/2}}\right)\frac{\al x}{2\lambda t^{\al/2+1}} $$
and 
$$ \textbf{1}_{[0,t_0]}(t)=\begin{cases} 1 & \text{ if }\,  0<t<t_0 \\
0 & \text{ else}.\end{cases} $$

\noindent Then function $v$ can be written as a convolution in the $t$--variable

$$ v(x,t)=K(x,t)\ast \textbf{1}_{[0,t]}. $$

\noindent This new way of expressing $v$ lead us to propose
 the following function

\begin{equation}{\label{c_2(x,t)}}
c_2(x,t)=K(x,t)\ast g(t)1_{[0,t]}=\int_0^t \MM_{\al/2}\left(\frac{x}{\lambda (t-\tau)^{\al/2}}\right)\frac{\al x}{2\lambda (t-\tau)^{\al/2+1}}g(\tau)d\tau
 \end{equation}
as a solution to problem (\ref{PVI-ICB}). \\

\noindent In order to prove this  assertion, let us enunciate the following lemma.

\begin{lemma}\label{derivacion alpha bajo el signo integral} Let $\KK(t-\tau)f(\tau)$ be a function that verifies the following conditions:

\begin{equation}\label{hip lema 1}
\hspace{-5.5 cm}\KK(t-\tau)f(\tau) \quad \text{ is a }\tau-\text{integrable function in } [0,t],
\end{equation}

\begin{equation}\label{hip lema 2}\hspace{-6.4 cm} \left|\frac{\partial }{\partial t}\KK(t-\tau)f(\tau)\right|\leq g(\tau) \, \in \, L^1[0,t], \quad \quad \quad \quad \quad
\end{equation}

\begin{equation}\label{hip lema 3} \hspace{-1.1 cm}\left|\frac{\partial }{\partial \eta}\KK(\eta-\tau)\frac{f(\tau)}{(t-\eta)^\al}\right| \, \in L^1(\Omega), \text{ where } \Omega=\{(\eta,\tau)\in \bbR^2 : \eta \in (0,t),\, 0\leq \tau \leq \eta \}
\end{equation}

\begin{equation}\label{hip lema 4} \hspace{-7.4 cm} 
\lim\limits_{\tau \nearrow \eta}\KK(\eta - \tau)f(\tau)=h(\eta) \, \in \, L^1(0,t).\quad \quad
\end{equation}

%\begin{figure}[h]
%\vspace{-4 cm}	
%\includegraphics[width=0.18\textwidth]{veamos_omega_este.eps}
%\end{figure}

%\vspace{0.3cm}
\noindent Then
 $$ _0D^{\al}_t \left( \int^{t}_{0} \KK(t-\tau)f(\tau)d\tau\right)=
 \int^{t}_{0}\left(_0D^{\al}_tK(t-\tau)\right)f(\tau)d\tau + _0I_t^{1-\al}\left( h(t) \right), $$

\noindent where $_0I^{1-\al}$ is the fractional integral of Riemann--Liouville of order $1-\al$ defined by 
$$_{0}I^{1-\alpha}h(t)=\frac{1}{\Gamma(1-\alpha)}\int^{t}_{0}(t-\tau)^{(1-\alpha)-1} h(\tau)d\tau. $$

\end{lemma}

%\begin{floatingfigure}[r]{30mm}
% \epsfysize=30mm
% \epsfbox{veamos_omega_este.eps}
% \end{floatingfigure}
%

\noindent \textit{Proof.}\\
 Due to (\ref{hip lema 2}) and (\ref{hip lema 4})

$$ \frac{d}{dt}\int^{t}_{0}\KK(t-\tau)f(\tau)d\tau= \int^{t}_{0}\frac{\partial }{\partial t}\KK(t-\tau)f(\tau)d\tau + \lim_{\tau \nearrow t}\KK(t-\tau)f(\tau)  . $$
\noindent Now

$$\hspace{-2cm}  _0D^{\al}_t \left( \int^{t}_{0} \KK(t-\tau)f(\tau)d\tau\right)=\frac{1}{\G(1-\al)}\int^{t}_{0}\frac{\frac{\partial}{\partial \eta}\int^{\eta}_{0}\KK(\eta - \tau)f(\tau)d\tau}{(t-\eta)^\al}d\eta= \quad \quad $$
 \begin{equation}\label{justif fubini lema 3} =\frac{1}{\G(1-\al)}\int^{t}_{0}\frac{1}{(t-\eta)^\al}\left[\int^{\eta}_{0}\frac{\partial}{\partial \eta}\KK(\eta - \tau)f(\tau)d\tau+\lim_{\tau \nearrow \eta}\KK(\eta-\tau)f(\tau)\right]d\eta .
 \end{equation}
Since (\ref{hip lema 3}) holds, (\ref{justif fubini lema 3}) is equal to
$$ \quad \int^{t}_{0}\frac{f(\tau)}{\G(1-\al)}\int^{t}_{\tau}\frac{\frac{\partial}{\partial \eta}K(\eta - \tau)}{(t-\eta)^\al}d\eta d\tau+ \frac{1}{\G(1-\al)}\int^{t}_{0}\frac{ \lim_{\tau \nearrow \eta}\KK(\eta-\tau)f(\tau)}{(t-\eta)^\al} d\eta .   $$

\noindent Substituting $s=\eta-\tau$,

$$ \frac{1}{\G(1-\al)} \int^{t}_{\tau}\frac{\frac{\partial}{\partial \eta}\KK(\eta - \tau)}{(t-\eta)^\al}d\eta=\frac{1}{\G(1-\al)} \int^{t-\tau}_{0}\frac{\frac{\partial}{\partial s}\KK(s)}{(t-\tau-s)^\al}ds=\,  _0D^{\al}_tK(t-\tau).$$

\noindent On the other hand,
$$ \frac{1}{\G(1-\al)}\int^{t}_{0}\frac{ \lim_{\tau \nearrow \eta}\KK(\eta-\tau)f(\tau)}{(t-\eta)^\al} d\eta=
\frac{1}{\G(1-\al)}\int^{t}_{0}\frac{h(\eta)}{(t-\eta)^\al} d\eta=\,_0 I_t^{1-\al}\left( h(t)\right)\, .$$

\noindent Hence
$$  _0D^{\al}_t \left( \int^{t}_{0} \KK(t-\tau)f(\tau)d\tau\right)=
 \int^{t}_{0}\left(\,_0D^{\al}_t\KK(t-\tau)\right)f(\tau)d\tau +\, _0I_t^{1-\al}\left( h(t)\right)  \, .$$
\qed

\noindent Now the purpose is to prove that the kernel  
\begin{equation}\label{kernel} \MM_{\al/2}\left(\frac{x}{\lambda (t-\tau)^{\al/2}}\right)\frac{\al x}{2 \lambda (t-\tau)^{\al/2+1}}g(\tau) \end{equation} verifies the hypothesis of Lemma \ref{derivacion alpha bajo el signo integral}.\\

\noindent $\bullet$  Hypothesis (\ref{hip lema 1}):  
	
\begin{equation}\label{hip 12-1} \hspace{0cm}\int_0^t \left|\MM_{\al/2}\left(\frac{x}{\lambda (t-\tau)^{\al/2}}\right)\frac{\al x}{2\lambda (t-\tau)^{\al/2+1}}g(\tau)\right|d\tau=\int_{x/\lambda t^{\al/2}}^\infty \MM_{\al/2}\left(y\right)\left|g\left(t-( \frac{x}{\lambda y})^{2/\al}\right)\right|dy  \end{equation}

\noindent We know that:\\

\begin{equation}\label{momento mainardi}
\hspace{-3cm}\int_0^\infty y^n \MM_{\al/2}(y)dy=\frac{\G(n+1)}{\G\left(\frac{\al}{2}n+1\right)} \quad \text{ for all } \al \in (0,2) \quad \quad \text{ (see \cite{FM-Analytical properties of the Wright function}).}
\end{equation}

\begin{equation}\label{M positiva}
\hspace{-1.5cm}\text{ The Mainardi function is a positive and decreasing function in } \bbR^+,  \text{(see \cite{RoSa1})}.
\end{equation}

\hspace{-0.3cm} $g$ is a bounded function in $[0,t]$, that is 

\begin{equation}\label{g<M}
\left|g(\tau)\right| \leq M \quad \forall \, \tau \in [0,t].
\end{equation}

\noindent Then (\ref{hip 12-1}) is convergent and (\ref{kernel}) is $\tau-integrable$ in $[0,t].$ 

\medskip
\noindent $\bullet$ Hypothesis (\ref{hip lema 2}):
$$\hspace{-7cm}\frac{\p}{\p t}\left[\MM_{\al/2}\left(\frac{x}{\lambda (t-\tau)^{\al/2}}\right)\frac{\al x}{2 \lambda (t-\tau)^{\al/2+1}}\right]g(\tau)=$$

$$\hspace{-5cm}=\WW\left(-\frac{x}{\lambda (t-\tau)^{\al/2}},-\frac{\al}{2},1-\al\right)\left(\frac{\al x}{2\lambda (t-\tau)^{\al/2+1}}\right)^2g(\tau)-$$
\begin{equation}\label{nucleo hip 2}
\hspace{5cm}-\MM_{\al/2}\left(\frac{x}{\lambda (t-\tau)^{\al/2}}\right)\frac{(\al/2+1)\al x}{2\lambda (t-\tau)^{\al/2+2}}g(\tau) \end{equation}
Applying Corollary \ref{acotacion en el inf W(-s,-al/2,1-al)}, there exists $\delta>0$ such that, for all $\tau \in (t-\delta, t)$, 
$$\hspace{-5cm}\left| \WW\left(-\frac{x}{\lambda (t-\tau)^{\al/2}},-\frac{\al}{2},1-\al\right)\left(\frac{\al x}{2\lambda (t-\tau)^{\al/2+1}}\right)^2g(\tau)\right|\leq $$
$$\hspace{2cm}\leq\left(c+d\frac{x}{\lambda (t-\tau)^{\al/2}}\right)\exp\left\{-b\left(\frac{x}{\lambda (t-\tau)^{\al/2}}\right)^{\frac{2}{2-\al}}\right\}\left(\frac{\al x}{2\lambda (t-\tau)^{\al/2+1}}\right)^2\left|g(\tau)\right|,$$
$$ \quad b>0, \quad c, \, d \, \text{constants}.$$
And this is an integrable function, in fact, making the substitution 
\begin{equation}\label{sustitucion}\frac{x}{\lambda (t-\tau)^{\al/2}}=r
\end{equation}
 and considering the inequality
\begin{equation}\label{desigualdad exp}
\exp\left\{-x\right\}\leq \frac{n!}{x^n}, \quad \quad \forall \quad  n \in \bbN, \quad \forall\,  x>0, 
\end{equation}
it results that

$$\hspace{-0.8cm} \int_{t-\delta}^t  \left(c+d\frac{x}{\lambda (t-\tau)^{\al/2}}\right)\exp\left\{-b\left(\frac{x}{\lambda (t-\tau)^{\al/2}}\right)^{2/2-\al}\right\}\left(\frac{\al x}{2\lambda (t-\tau)^{\al/2+1}}\right)^2\left|g(\tau)\right|d\tau= $$
$$\hspace{-3.5cm} =\int_{\frac{x}{\lambda \delta^{\al/2}}}^\infty\left(c+dr\right)\exp\left\{-br^{2/2-\al}\right\}\frac{\al}{2} r\left(\frac{r\lambda}{x}\right)^{2/\al}\left|g\left(t-\left( \frac{x}{\lambda r} \right)^{2/\al}\right)\right|dr \leq $$
\begin{equation}\label{int-imp} \leq C_{x,\la,\al}M\int_{\frac{x}{\lambda \delta^{\al/2}}}^\infty \left(c+dr\right) r^{1+2/\al} \exp\left\{-br^{2/2-\al}\right\}dr \leq C_{x,\la,\al}M\int_{\frac{x}{\lambda \delta^{\al/2}}}^\infty \frac{\left(c+dr\right) r^{1+2/\al}}{b^n r^{2n/2-\al}} n! \, dr.\end{equation} 

\noindent It it easy to see that for any $\al \in (0,1)$, there exists $n\in \bbN$ such that (\ref{int-imp}) is convergent. For example,  if $\al= 1/4$ we can take $n= 10 $. \\
Then, the first term of the sum (\ref{nucleo hip 2}) is bounded by an integrable function. \\

\noindent Let us consider the second term of sum (\ref{nucleo hip 2}). Making the substitution (\ref{sustitucion}) and taking into account that the Mainardi function is a positive function, we have 
 
$$\hspace{-6.5cm} \int_0^t\left|\MM_{\al/2}\left(\frac{x}{\lambda (t-\tau)^{\al/2}}\right)\frac{(\al/2+1)\al x}{2\lambda (t-\tau)^{\al/2+2}} g(\tau)\right|d\tau\leq$$
\begin{equation}\hspace{1cm} \leq M\int_{x/\lambda t^{\al/2}}^{\infty} (\al/2+1)\MM_{\al/2}\left( r \right)\left(\frac{\lambda r}{x}\right)^{2/\al}dr \leq  MC_{x,\lambda,\al}\int_0^\infty \MM_{\al/2}(r) r^{2/\al}dr.  \end{equation}    

\noindent Now, for any $\al \in (0,1)$, there exists $k \in \bbN$ such that $ \frac{1}{k}<\frac{\al}{2}<1 .$
Then, using (\ref{momento mainardi}), it yields 
\begin{equation}\label{acot segunda parte}
\hspace{0cm} MC_{x,\lambda,\al}\int_0^\infty \MM_{\al/2}(r) r^{2/\al}dr \leq  MC_{x,\lambda,\al}\int_0^\infty \MM_{\al/2}(r) r^{k}dr=MC_{x,\lambda,\al} \frac{\G(k+1)}{\G\left(\frac{\al}{2}k +1\right)}.   \end{equation}

\noindent $\bullet$ Hypothesis (\ref{hip lema 3}): \\
We have to prove that 
$$ \frac{\p}{\p \eta}\left[\MM_{\al/2}\left(\frac{x}{\lambda (\eta-\tau)^{\al/2}}\right)\frac{\al x}{2\lambda (\eta-\tau)^{\al/2+1}}\right]\frac{g(\tau)}{(t-\eta)^\al}  \, \in L^1(\Omega), $$
 where  $\Omega=\{(\eta,\tau)\in \bbR^2 : \, \eta \in (0,t),\, 0\leq \tau \leq \eta \}$. Or equivalently,
$$\hspace{-5.5cm} \frac{\p}{\p \eta}\left[\MM_{\al/2}\left(\frac{x}{\lambda (\eta-\tau)^{\al/2}}\right)\right]\left(\frac{\al x}{2\lambda (\eta-\tau)^{\al/2+1}}\right)\frac{g(\tau)}{(t-\eta)^\al}- $$
$$\hspace{4cm} - \MM_{\al/2}\left(\frac{x}{\lambda (\eta-\tau)^{\al/2}}\right)\frac{(\al/2+1)\al x}{2\lambda(\eta-\tau)^{\al/2+2}}\frac{g(\tau)}{(t-\eta)^\al}=I+II  \, \in L^1(\Omega). $$

\noindent Reasoning like in the previous item, using Corollary \ref{acotacion en el inf W(-s,-al/2,1-al)}, inequality (\ref{desigualdad exp}), Corollary \ref{acotacion M} and Tonelli's theorem (see \cite{Brezis}, p. 55), the following assertions are true:
\begin{equation}\label{hip Tonelli 1} \int_0^\eta\left|I + II\right|d\tau <\infty \, \forall\, \eta \in (0,t). 
\end{equation}  
Taking $\delta$ small according to Corollary \ref{acotacion en el inf W(-s,-al/2,1-al)}, 

\begin{equation}\label{hip Tonelli 2-1} \int_0^t\int_0^\eta\left|I + II\right|d\tau d\eta = \int_0^t\int_0^{\eta-\delta} \left|I + II\right|d\tau d\eta + \int_0^t\int_{\eta-\delta}^\eta\left|I + II\right|d\tau d\eta. \end{equation}  

\noindent Now, noting that
\begin{equation}
\frac{\p}{\p \eta}\left[\MM_{\al/2}\left(\frac{x}{\lambda (\eta-\tau)^{\al/2}}\right)\right]=-\frac{\p}{\p \tau}\left[\MM_{\al/2}\left(\frac{x}{\lambda (\eta-\tau)^{\al/2}}\right)\right]
\end{equation}
and that
\begin{equation}
-\frac{\p}{\p \tau}\MM_{\al/2}\left(\frac{x}{\lambda (\eta-\tau)^{\al/2}}\right)>0 \quad \text{(it is a consequence of Lemma 4.2 from \cite{RoSa1}).}
\end{equation}

\noindent Let be $M$ defined in $(\ref{g<M})$ and $C$ any constant depending on $\delta, \al, x$ or $n$. Then

$$%\hspace{-9.5cm}
\int_0^t\int_0^{\eta-\delta}\left|I\right|d\tau d\eta \leq \int_0^t \frac{M\,C}{(t-\eta)^\al \delta^{\al/2+1}}\int_0^{\eta -\delta} -\frac{\p}{\p \tau}\MM_{\al/2}\left(\frac{x}{\lambda (\eta-\tau)^{\al/2}}\right)d\tau d\eta = $$
%$$\hspace{3cm}\leq \int_0^t \frac{M}{(t-\eta)^\al}\int_0^{\eta -\delta}\left| 
%\frac{\p}{\p \eta}\left[\MM_{\al/2}\left(\frac{x}{\lambda (\eta-\tau)^{\al/2}}\right)\right]\left(\frac{\al x}{2\lambda (\eta-\tau)^{\al/2+1}}\right)\right|d\tau d\eta= $$  
%
%$$ \leq \int_0^t \frac{M}{(t-\eta)^\al}\int_0^{\eta -\delta}\left|-
%\frac{\p}{\p \tau}\MM_{\al/2}\left(\frac{x}{\lambda (\eta-\tau)^{\al/2}}\right)\right|\left(\frac{\al x}{2\lambda (\eta-\tau)^{\al/2+1}}\right)d\tau d\eta\leq $$

%$$ =\int_0^t \frac{M\,C}{(t-\eta)^\al \delta^{\al/2+1}}\left[ - \left. \MM_{\al/2}\left(\frac{x}{\lambda (\eta-\tau)^{\al/2}}\right)\right|_0^{\eta -\delta}\right] d\eta = $$
$$\hspace{3cm}=\int_0^t \frac{M\,C}{(t-\eta)^\al \delta^{\al/2+1}}\left[  \MM_{\al/2}\left(\frac{x}{\lambda \eta^{\al/2}}\right)- \MM_{\al/2}\left(\frac{x}{\lambda \delta^{\al/2}}\right) \right] d\eta \leq $$
$$\hspace{-0.5cm}\leq 
\frac{2M\,C}{\G\left(1-\frac{\al}{2}\right)\delta^{\al/2+1}}\int_0^t \frac{1}{(t-\eta)^\al } d\eta <\infty$$
due to Lemma 4.2 \cite{RoSa1} and that $\al \in (0,1)$. 
Then, 
\begin{equation}\label{Tonelli 2-1}
\int_0^t\int_0^{\eta-\delta}\left|I\right|d\tau d\eta<\infty. 
\end{equation}
On the other hand, 
$$\hspace{-0.8cm}
\int_0^t\int_0^{\eta-\delta}\left|II\right|d\tau d\eta =  \int_0^t \int_{0}^{\eta-\delta}\left| 
\MM_{\al/2}\left(\frac{x}{\lambda (\eta-\tau)^{\al/2}}\right)\frac{(\al/2+1)\al x}{2\lambda(t-\tau)^{\al/2+2}}\frac{g(\tau)}{(t-\eta)^\al}\right|d\tau d\eta= $$  
$$ \hspace{2.8cm}= \int_0^t \frac{1}{(t-\eta)^\al}\int_{0}^{\eta-\delta}\left| 
\MM_{\al/2}\left(\frac{x}{\lambda (\eta-\tau)^{\al/2}}\right)\frac{(\al/2+1)\al x}{2\lambda(t-\tau)^{\al/2+2}}g(\tau)\right|d\tau d\eta. $$
We proved in the previous item that
$$ \int_0^\eta\left|\MM_{\al/2}\left(\frac{x}{\lambda (\eta-\tau)^{\al/2}}\right)\frac{(\al/2+1)\al x}{2\lambda (\eta-\tau)^{\al/2+2}} g(\tau)\right|d\tau <\infty.  $$

\noindent Recalling that $\al \in (0,1)$, it yields that 

\begin{equation}\label{Tonelli 2-2}
\int_0^t\int_0^{\eta-\delta}\left|II\right|d\tau d\eta <\infty.
\end{equation}  
Then
\begin{equation}\label{hip-Tonelli 2a}
\int_0^t\int_0^{\eta-\delta}\left|I + II\right|d\tau d\eta <\infty.
 \end{equation}

\noindent Proceeding like in item $2$, it can be proved that 
\begin{equation}\label{hip-Tonelli 2b} \int_0^t\int_{\eta-\delta}^{\eta}\left|II\right|d\tau d\eta <\infty. 
\end{equation}

\noindent Finally, (\ref{hip-Tonelli 2a}) and (\ref{hip-Tonelli 2b}) yield 
\begin{equation}\label{hip-Tonelli 2 hecha}
\int_0^t\int_0^\eta\left|I + II\right|d\tau d\eta <\infty.
 \end{equation}

\noindent From $(\ref{Tonelli 2-1})$ and $(\ref{hip-Tonelli 2 hecha})$, Tonelli's theorem holds and
$$ \frac{\p}{\p \eta}\left[\MM_{\al/2}\left(\frac{x}{\lambda (\eta-\tau)^{\al/2}}\right)\frac{x}{\lambda (\eta-\tau)^{\al/2+1}}\frac{\al}{2}\right]\frac{g(\tau)}{(t-\eta)^\al}  \, \in L^1(\Omega). $$

\noindent $\bullet$ Hypothesis (\ref{hip lema 4}):\\

\noindent Let us prove that \begin{equation}\label{lim 1} \lim\limits_{\tau \nearrow \eta}\MM_{\al/2}\left(\frac{x}{\lambda (\eta-\tau)^{\al/2}}\right)\frac{\al x}{2 \lambda (\eta-\tau)^{\al/2+1}}g(\tau)=0.\end{equation} 
Note that, due to  (\ref{momento mainardi}) and  (\ref{g<M}),
$$\hspace{-7cm} \lim\limits_{\tau \nearrow \eta}\MM_{\al/2}\left(\frac{x}{\lambda (\eta-\tau)^{\al/2}}\right)\frac{\al x}{2 \lambda (\eta-\tau)^{\al/2+1}}g(\tau)= $$
$$ =\lim\limits_{s\rightarrow^+ 0}\MM_{\al/2}\left(\frac{x}{\lambda s^{\al/2}}\right)\frac{\al x}{2 \lambda s^{\al/2+1}}g(\eta-s)=\lim\limits_{y\rightarrow \infty}\MM_{\al/2}\left(\frac{x}{\lambda }y^{\al/2}\right)\frac{\al x y^{\al/2+1}}{2 \lambda }g\left(\eta-\frac{1}{y^{\al/2+1}}\right)=0 .$$

\bigskip
\noindent Finally, we can apply Lemma \ref{derivacion alpha bajo el signo integral} to kernel (\ref{c_2(x,t)}). Then,\\

\newpage
$$\hspace{-2.5cm}  _0D^{\al} c_2(x,t)=\,_0D^{\al}\left( \int_0^t \MM_{\al/2}\left(\frac{x}{\lambda (t-\tau)^{\al/2}}\right)\frac{\al x}{2\lambda (t-\tau)^{\al/2+1}}g(\tau)d\tau \right)=\quad
$$
$$\hspace{1cm}  = \int^{t}_{0}\, _0D^{\al} \left( \MM_{\al/2}\left(\frac{x}{\lambda (t-\tau)^{\al/2}}\right)\frac{\al x}{2\lambda (t-\tau)^{\al/2+1}}g(\tau)\right)d\tau  +\, I_t^{1-\al}\left( 0 \right) =
  $$
$$ \hspace{-2cm}  =\int^{t}_{0} \,_0D^{\al} \left(-\frac{\p}{\p t}\WW\left(\frac{-x}{\lambda (t-\tau)^{\al/2}},-\frac{\al}{2},1\right)\right)g(\tau)d\tau=$$	
	\begin{equation}\label{aplic lema-1}
	\hspace{-2.5cm} =	\int^{t}_{0}\, _0D^{\al+1} \left(-\WW\left(\frac{-x}{\lambda (t-\tau)^{\al/2}},-\frac{\al}{2},1\right)\right)g(\tau)d\tau.\end{equation}
	
	\noindent From Corollary \ref{lim M(x)=0}, $\frac{\p }{\p \tau}\WW\left(\frac{-x}{\lambda 0^{+}},-\frac{\al}{2},1\right)=0$. Then Caputo derivative  commutes  (see eq. (2.145) of \cite{Podlubny}) and 
	$$\hspace{-1.8cm} _0D^{\al+1} \left(-\WW\left(\frac{-x}{\lambda (t-\tau)^{\al/2}},-\frac{\al}{2},1\right)\right)=\,_0D^{1+\al} \left(-\WW\left(\frac{-x}{\lambda (t-\tau)^{\al/2}},-\frac{\al}{2},1\right)\right)=$$
	$$\hspace{4.4cm}=-\frac{\p}{\p t}\left(\,_0D^{\al}\WW\left(\frac{-x}{\lambda (t-\tau)^{\al/2}},-\frac{\al}{2},1\right)\right). $$
	
\noindent Now, recalling that $\WW\left(\frac{-x}{\lambda (t-\tau)^{\al/2}},-\frac{\al}{2},1\right)$ is a solution of problem (\ref{PVI1}), it results that 

$$\hspace{-2.2cm}-\frac{\p}{\p t}\left(\,_0D^{\al} \WW\left(\frac{-x}{\lambda (t-\tau)^{\al/2}},-\frac{\al}{2},1\right)\right)= -\frac{\p}{\p t} \frac{\p}{\p x^2} \WW\left(\frac{-x}{\lambda (t-\tau)^{\al/2}},-\frac{\al}{2},1\right)= $$
\begin{equation}\label{aplic lema-2}\hspace{-0.1cm}=\frac{\p}{\p x^2}\left(-\frac{\p}{\p t}  \WW\left(\frac{-x}{\lambda (t-\tau)^{\al/2}},-\frac{\al}{2},1\right)\right)=\frac{\p}{\p x^2}\MM_{\al/2}\left(\frac{x}{\lambda (t-\tau)^{\al/2}}\right)\frac{x}{\lambda (t-\tau)^{\al/2+1}}\frac{\al}{2}\end{equation}

\noindent Replacing (\ref{aplic lema-2}) in (\ref{aplic lema-1}), we have 

\begin{equation}
 \hspace{-3.5 cm}_0D^{\al} c_2(x,t)=
	 \int^{t}_{0}\frac{\p}{\p x^2}\MM_{\al/2}\left(\frac{x}{\lambda (t-\tau)^{\al/2}}\right)\frac{x}{\lambda (t-\tau)^{\al/2+1}}\frac{\al}{2}g(\tau)d\tau=$$
	$$\hspace{0.5cm}=\frac{\p}{\p x^2}\int^{t}_{0}\MM_{\al/2}\left(\frac{x}{\lambda (t-\tau)^{\al/2}}\right)\frac{x}{\lambda (t-\tau)^{\al/2+1}}\frac{\al}{2}g(\tau)d\tau=\frac{\p}{\p x^2}c_2(x,t).
\end{equation}
Therefore $c_2(x,t)$ verifies the fractional diffusion equation.

\bigskip
%
%Veamos que $c_2(x,t)$ verifica las condiciones de borde. \\
%Pero antes, enunciamos la siguiente proposici\'on que nos ser\'a de gran utilidad. 

\begin{proposition} \label{salto derecho} The following limits holds:
%Si $g\equiv 1$ en (\ref{c_2(x,t)}), entonces 
\begin{enumerate}
	\item $\lim\limits_{x\searrow 0}\int^{t}_{0}\MM_{\al/2}\left(\frac{x}{\lambda (t-\tau)^{\al/2}}\right)\frac{\al x}{2\lambda (t-\tau)^{\al/2+1}}d\tau=1 .$
	\item $\lim\limits_{t\searrow 0}\int^{t}_{0}\MM_{\al/2}\left(\frac{x}{\lambda (t-\tau)^{\al/2}}\right)\frac{\al x}{2 \lambda (t-\tau)^{\al/2+1}}d\tau=0 .$
\end{enumerate}  
\end{proposition}
\noindent \textit{Proof.}\\ 
1. Taking $n=0$ in (\ref{momento mainardi}), it yields $\int_0^\infty \MM_{\al/2}(u)du=1$. \\
Making the substitution   (\ref{sustitucion}) and applying Corollary \ref{acotacion M}, we have
$$ \left| \int^{t}_{0}\MM_{\al/2}\left(\frac{x}{\lambda (t-\tau)^{\al/2}}\right)\frac{\al x}{2\lambda (t-\tau)^{\al/2+1}}d\tau -1 \right| =\left| \int_{x/\lambda t^{\al/2}}^\infty \MM_{\al/2}(u) du - \int_0^\infty \MM_{\al/2}(u)du \right|\leq  $$
$$ \hspace{4cm}\leq \int_0^{x/\lambda t^{\al/2}} \MM_{\al/2}(u) du \leq \frac{1}{\G\left(1-\frac{\al}{2}\right)}\frac{x}{\lambda t^{\al/2}} \rightarrow 0 \quad \text{ if } \quad x\searrow 0.$$
\noindent 2. It is a consequence of applying substitution (\ref{sustitucion}). 
\qed

\medskip
\noindent Let us check now the border conditions.\\

\noindent $\bullet$ $c_2(x,0)=\lim\limits_{t\searrow 0}c_2(x,t)$.\\
		From Proposition \ref{salto derecho} and (\ref{g<M}), 
	$$\hspace{-5cm} 0\leq \lim\limits_{t\searrow 0}\left|\int^{t}_{0}\MM_{\al/2}\left(\frac{x}{\lambda (t-\tau)^{\al/2}}\right)\frac{\al x}{2 \lambda (t-\tau)^{\al/2+1}}g(\tau)d\tau\right|\leq $$
	$$\hspace{-4.2cm}  \leq 	\lim\limits_{t\searrow 0}\int^{t}_{0}\MM_{\al/2}\left(\frac{x}{\lambda (t-\tau)^{\al/2}}\right)\frac{\al x}{2 \lambda (t-\tau)^{\al/2+1}}\frac{\al}{2}\left|g(\tau)\right|d\tau \leq$$
	$$\hspace{-4.3cm}  \leq \lim\limits_{t\searrow 0}M\int^{t}_{0}\MM_{\al/2}\left(\frac{x}{\lambda (t-\tau)^{\al/2}}\right)\frac{\al x}{2 \lambda (t-\tau)^{\al/2+1}}\frac{\al}{2}d\tau=0 $$
	\begin{equation} \label{cond de borde PVIC t-0 }
	\therefore \quad c_2(x,0)=0.
	\end{equation}

\noindent $\bullet$	$ c_2(0,t)=\lim\limits_{x\searrow 0}c_2(x,t)$.\\
\noindent 	Note that
 	
$$\hspace{-6.8cm} \int^{t}_{0}\MM_{\al/2}\left(\frac{\al x}{2\lambda (t-\tau)^{\al/2}}\right)\frac{x}{\lambda (t-\tau)^{\al/2+1}}\frac{\al}{2}g(\tau)d\tau=$$
$$ \hspace{-3.8cm} =\int^{t}_{0}\MM_{\al/2}\left(\frac{x}{\lambda (t-\tau)^{\al/2}}\right)\frac{\al x}{2\lambda (t-\tau)^{\al/2+1}}\frac{\al}{2}\left[g(\tau)-g(t)+g(t)\right]d\tau= $$
$$\hspace{-5.3cm} =\int^{t}_{0}\MM_{\al/2}\left(\frac{x}{\lambda (t-\tau)^{\al/2}}\right)\frac{\al x}{2\lambda (t-\tau)^{\al/2+1}}\left[g(\tau)-g(t)\right]d\tau+$$
$$\hspace{6.5cm} +\int^{t}_{0}\MM_{\al/2}\left(\frac{x}{\lambda (t-\tau)^{\al/2}}\right)\frac{\al x}{2\lambda (t-\tau)^{\al/2+1}}g(t)d\tau .$$

\noindent Applying Proposition \ref{salto derecho} to the second member of the sum, 
\newpage
$$\hspace{-6cm} \lim\limits_{x\searrow 0}\int^{t}_{0}\MM_{\al/2}\left(\frac{x}{\lambda (t-\tau)^{\al/2}}\right)\frac{\al x}{2\lambda (t-\tau)^{\al/2+1}}g(t)d\tau =$$
$$\hspace{4.5cm}=g(t)\lim\limits_{x\searrow 0}\int^{t}_{0}\MM_{\al/2}\left(\frac{x}{\lambda (t-\tau)^{\al/2}}\right)\frac{\al x}{2\lambda (t-\tau)^{\al/2+1}}d\tau=g(t). $$

\noindent Let be \begin{equation}\label{I}I=\int^{t}_{0}\MM_{\al/2}\left(\frac{x}{\lambda (t-\tau)^{\al/2}}\right)\frac{x}{\lambda (t-\tau)^{\al/2+1}}\frac{\al}{2}\left[g(\tau)-g(t)\right]d\tau.\end{equation}
The next goal is to prove that $\lim\limits_{x\searrow 0}I=0$. \\

\noindent Let  $\delta>0 $ to be determined. 
 $$\hspace{-4.3cm} I=\int^{t-\delta}_{0}\MM_{\al/2}\left(\frac{x}{\lambda (t-\tau)^{\al/2}}\right)\frac{\al x}{2 \lambda (t-\tau)^{\al/2+1}}\left[g(\tau)-g(t)\right]d\tau + $$
\begin{equation}\label{}\hspace{2 cm} + \int_{t-\delta}^{t}\MM_{\al/2}\left(\frac{x}{\lambda (t-\tau)^{\al/2}}\right)\frac{\al x}{2 \lambda (t-\tau)^{\al/2+1}}\frac{\al}{2}\left[g(\tau)-g(t)\right]d\tau= I_1+I_2 \end{equation}

\noindent Applying Corollary \ref{acotacion M},
$$\hspace{-3.8cm} |I_1|=\left|  \int^{t-\delta}_{0}\MM_{\al/2}\left(\frac{x}{\lambda (t-\tau)^{\al/2}}\right)\frac{\al x}{2\lambda (t-\tau)^{\al/2+1}}\left[g(\tau)-g(t)\right]d\tau  \right|\leq  $$
$$ \leq  \int^{t-\delta}_{0}\frac{1}{\G\left(1-\frac{\al}{2}\right)}\frac{\al x}{2 \lambda (t-\tau)^{\al/2+1}}\frac{\al}{2}|g(\tau)-g(t)|d\tau\leq \frac{1}{\G\left(1-\frac{\al}{2}\right)}\frac{\al x}{2 \lambda \delta^{\al/2+1}}\int^{t-\delta}_{0}|g(\tau)-g(t)|d\tau  $$
\begin{equation}\label{I1<e}\hspace{-1.7cm} \leq  \left(\frac{1}{\G\left(1-\frac{\al}{2}\right)}\frac{\al}{2\lambda \delta^{\al/2+1}}\int^{t}_{0}|g(\tau)-g(t)|d\tau \right)x=C_{t,\delta,\al}x <\frac{\epsilon}{2},  \quad \text{if } \quad x<\frac{\e}{2C_{t,\delta,\al}}.\end{equation}

\noindent Now,  $g$ is  continuous in $t$. Then, for $\e>0$ given, there exists $\delta>0$ such that \\
$|g(\tau)-g(t)|<\frac{\e}{2}$ if $|t-\tau|<\delta$. Using this fact and making the substitution (\ref{sustitucion}), 

\begin{equation}\label{I2<e} |I_2|\leq  \frac{\e}{2}\int_{t-\delta}^{t}\MM_{\al/2}\left(\frac{x}{\lambda (t-\tau)^{\al/2}}\right)\frac{x}{\lambda (t-\tau)^{\al/2+1}}\frac{\al}{2}d\tau = \frac{\e}{2} \int_{x/\lambda \delta^{\al/2}}^\infty \MM_{\al/2}(u) du < \frac{\e}{2}. \end{equation}

\noindent From (\ref{I1<e}) and (\ref{I2<e}), it results that $ |I|<\e $, for every $\e>0$.  Then, 
$$ \lim\limits_{x\searrow 0}\int^{t}_{0}\MM_{\al/2}\left(\frac{x}{\lambda (t-\tau)^{\al/2}}\right)\frac{x}{\lambda (t-\tau)^{\al/2+1}}\frac{\al}{2}\left[g(\tau)-g(t)\right]d\tau =0,$$

\noindent and we can assure that

\begin{equation} \label{cond de borde PVIC x-0} 
 c_2(0,t)=g(t).
\end{equation}

\begin{theorem}\label{sol PVIC-IC} Let be $f $ a continuous bounded function in $\bbR^+_0$ and  $g$ a continuous function in $[0,T)$ 
. Then 

$$c(x,t)=\int^{\infty}_{0}\frac{1}{2 \lambda t^{\al/2}}\left[ \MM_{\al/2}\left(\frac{|x-\xi|}{\lambda t^{\al/2}}\right) - \MM_{\al/2}\left(\frac{|x+\xi|}{\lambda t^{\al/2}}\right) \right]f(\xi) d\xi+\quad \quad \quad \quad \quad \quad \quad \quad$$
\begin{equation}\label{c(x,t) solución}
\quad \quad \quad \quad\quad \quad \quad\quad \quad \quad \quad\quad \quad \quad  \int_0^t \MM_{\al/2}\left(\frac{x}{\lambda (t-\tau)^{\al/2}}\right)\frac{\al x}{2 \lambda (t-\tau)^{\al/2+1}}g(\tau)d\tau
\end{equation}

\noindent  is a solution to problem

\begin{equation}{\label{PVI-IC-teo}}
\left\{\begin{array}{lll}
          D^{\al} c(x,t)=\lambda^2\dfrac{\partial^2c }{\partial x^2}(x,t) &  & 0<x<\infty, \, 0<t<T, \, 0<\al<1, \\

          c(x,0)=f(x) & \quad  \quad  & 0<x<\infty, \\

          c(0,t)=g(t) & \quad  \quad  & 0<t<T. \\

       \end{array}\right.  \end{equation}
\end{theorem}

\begin{theorem}The limit when $\al\nearrow 1$ of the solution to problem
\begin{equation}\label{PVIC-IC-fnula}\left\{\begin{array}{lll}
          D^{\al} c_2(x,t)=\lambda^2\dfrac{\partial^2c_2 }{\partial x^2}(x,t) &  & 0<x<\infty, \, 0<t<T, \, 0<\al<1, \\
          c_2(x,0)=0 & \quad  \quad  & 0<x<\infty, \\

          c_2(0,t)=g(t) & \quad  \quad  & 0<t<T, \\

                                       \end{array}\right. %\end{equation}
\end{equation} is the classical solution to the analogous  problem when $\al=1$ and we recover the heat equation
\begin{equation}{\label{PVI-ec del calor}}
\left\{\begin{array}{lll}
          \frac{\partial}{\partial t}w(x,t)=\lambda^2\dfrac{\partial^2c }{\partial x^2}(x,t) &  & 0<x<\infty, \, 0<t<T, \\
          w(x,0)=0 & \quad  \quad  & 0<x<\infty, \\

          w(0,t)=g(t) & \quad  \quad  & 0<t<T. \\

                                       \end{array}\right.\end{equation}
\end{theorem}
\proof
Let be 
\begin{equation}\label{sol PVI-IC ec del calor}
w(x,t)=\int^{t}_{0}\frac{e^{-\frac{x^2}{4\lambda^2(t-\tau)}}}{2\sqrt{\pi}}
\frac{x}{\lambda(t-\tau)^{3/2}}g(\tau)d\tau 
\end{equation}
the solution of problem (\ref{PVI-ec del calor}) (see \cite{Cannon}) and let be $c^{\al}_2$ the solution of problem (\ref{PVIC-IC-fnula}) given by Theorem \ref{sol PVIC-IC}, 
\begin{equation}\label{sol PVIC-IC-nula} c^{\al}_2(x,t)= \int_0^t \MM_{\al/2}\left(\frac{x}{\lambda (t-\tau)^{\al/2}}\right)\frac{x}{\lambda (t-\tau)^{\al/2+1}}\frac{\al}{2}g(\tau)d\tau. \end{equation}
\noindent Applying Lebesgue Convergence Theorem and Lemma \ref{conv M al/2 cuando al tiende a 1}, 
\newpage
$$\lim_{\al \nearrow 1 } c^{\al}_2(x,t) =\lim_{\al \nearrow 1 }\left\{ \int_0^t \MM_{\al/2}\left(\frac{x}{\lambda (t-\tau)^{\al/2}}\right)\frac{\al x}{2\lambda (t-\tau)^{\al/2+1}}g(\tau)d\tau \right\}=$$
$$\hspace{1.5cm}=\int_0^t \lim_{\al \nearrow 1 } \MM_{\al/2}\left(\frac{x}{\lambda (t-\tau)^{\al/2}}\right)\frac{\al x}{2\lambda (t-\tau)^{\al/2+1}}g(\tau)d\tau=$$
$$ \hspace{-0.8cm}= \int^{t}_{0}\frac{e^{-\frac{x^2}{4(t-\tau)}}}{2 \sqrt{ \pi}}\frac{x}{\lambda (t-\tau)^{3/2}}g(\tau)d\tau = w(x,t). $$

\endproof

\section{Solving the Initial-Boundary-Value Problem for the  Time-Fractional Diffusion Equation in the Quarter Plane with null Flux- Boundary Condition.}
\noindent In this section,  problem

%We consider now the following problem
\begin{equation}\label{p1}\left\{
\begin{array}{lll}
_0D^{\al} c_3(x,t)=\lambda^2\dfrac{\partial^2c_3 }{\partial x^2}(x,t) & & 0<x<\infty, \, 0<t<T, \, 0<\al<1 \\
          \dfrac{\partial c_3 }{\partial x}(0,t)=0  & &     0<t<T \\
           c_3(x,0)=f(x)   & & 0<x<\infty \\
           %\displaystyle\lim_{x\rightarrow +\infty}c(x,t)=0 & & t>0
\end{array}\right.
\end{equation}

%\noindent will be studied, where $f$....???????????\\

\noindent Let us consider the following auxiliary problem

$$
\left\{
\begin{array}{lll}
_0D^{\al} z(x,t)=\lambda^2\dfrac{\partial^2c }{\partial x^2}(x,t) & &  -\infty <x<\infty, \, 0<t<T, \, 0<\al<1 \\
           z(x,0)=\tilde{f} (x)   & &  -\infty <x<\infty  \\
           %\displaystyle\lim_{x\rightarrow \pm\infty}c(x,t)=0 & & t>0
\end{array}\right.
$$
where  $\tilde{f}$ is an even extension of $f$.\\
\noindent This problem was solved in \cite{FM2:1994} and its solution is given by the following function
$$
c_3(x,t)=\int^{\infty}_{-\infty} \frac{t^{-\frac{\al}{2}}}{2\lambda}
\MM_{\al/2}\left(\frac{|x-\xi |}{
\lambda t^{\frac{\al}{2}}}\right) \tilde{f}(\xi )d\xi= \hspace{8 cm} $$
\begin{equation}{\label{auxc}} \hspace{4cm}= \frac{1}{2\lambda t^{\frac{\al}{2}}}\int^{\infty}_{0}\left[\MM_{\al/2}\left(\frac{x+\xi}{\lambda t^{\frac{\al}{2}}}\right) + \MM_{\al/2}\left(\frac{|x-\xi |}{\lambda t^{\frac{\al}{2}}}\right)\right]
 f(\xi )d\xi . \end{equation}

\noindent The next goal is to prove that  
\begin{equation}\label{cond de flujo nula}
\displaystyle\lim_{x\searrow 0}
\frac{\partial}{\partial x}c_3(x,t)=0.
\end{equation}
 $$\dfrac{\partial c_3}{\partial x}(0,t)=\displaystyle\lim_{x\searrow 0}
\frac{1}{2\lambda t^{\frac{\al}{2}}}\frac{\partial}{\partial x}\int^{\infty}_{0}\left[\MM_{\al/2}\left(\frac{x+\xi}{
\lambda t^{\frac{\al}{2}}}\right) + \MM_{\al/2}\left(\frac{|x-\xi |}{\lambda t^{\frac{\al}{2}}}\right)\right]
 f(\xi )d\xi=\quad \quad  \quad \quad \quad \quad  \quad \quad \quad \quad  \quad \quad
$$
$$ \hspace{-1.3cm}=\displaystyle\lim_{x\searrow 0}
\frac{1}{2\lambda t^{\frac{\al}{2}}}\frac{\partial}{\partial x}\left\{\int^{x}_{0}\left[\MM_{\al/2}\left(\frac{x+\xi}{
\lambda t^{\frac{\al}{2}}}\right)+ \MM_{\al/2}\left(\frac{x-\xi}{
\lambda t^{\frac{\al}{2}}}\right)\right] f(\xi )d\xi
 + \right.$$
\begin{equation}{\label{dv/dx(0,t)}}
 \hspace{3.5cm} \left.+\int^{\infty}_{x}\left[\MM_{\al/2}\left(\frac{x+\xi}{\lambda t^{\frac{\al}{2}}}\right) +\MM_{\al/2}\left(\frac{\xi -x}{\lambda t^{\frac{\al}{2}}}\right)\right] f(\xi )d\xi \right\} .
\quad \quad   \quad \quad \end{equation}

\bigskip

\noindent Due to the continuity of the Mainardi function, Corollary \ref{acotacion en el inf W(-s,-al/2,1-al)} and $f \in \CC(\bbR^+)\cap L^\infty(\bbR^+)$,  the next euqalities are true:    

$$ \hspace{-6.5cm}\frac{\partial}{\partial x}\int^{x}_{0}\left[M_{\frac{\al}{2}}\left(\frac{x+\xi}{
\lambda t^{\frac{\al}{2}}}\right) + \MM_{\al/2}\left(\frac{x-\xi}{
\lambda t^{\frac{\al}{2}}}\right)\right] f(\xi )d\xi=
$$

$$\hspace{-5.8cm}=\int^{x}_{0}\frac{\partial}{\partial x} \left[\MM_{\al/2}\left(\frac{x+\xi}{
\lambda t^{\frac{\al}{2}}}\right) + \MM_{\al/2}\left(\frac{x-\xi}{
\lambda t^{\frac{\al}{2}}}\right)\right] f(\xi )d\xi  +$$
\begin{equation}{\label{d/dx(int0_x)}}\hspace{5cm} +\displaystyle\lim_{\xi \nearrow x}
\left[\MM_{\al/2}\left(\frac{x+\xi}{
\lambda t^{\frac{\al}{2}}}\right) + \MM_{\al/2}\left(\frac{x-\xi}{
\lambda t^{\frac{\al}{2}}}\right)\right]f(\xi )
\end{equation}

\noindent and

$$\hspace{-6cm} \frac{\partial}{\partial x} \int^{\infty}_{x}\left[\MM_{\al/2}\left(\frac{x+\xi}{
\lambda t^{\frac{\al}{2}}}\right) +\MM_{\al/2}\left(\frac{\xi-x}{
\lambda t^{\frac{\al}{2}}}\right)\right] f(\xi )d\xi=
 $$

$$ \hspace{-5.5cm}=-\int^{x}_{\infty}\frac{\partial}{\partial x}
\left[\MM_{\al/2}\left(\frac{x+\xi}{
\lambda t^{\frac{\al}{2}}}\right) +\MM_{\al/2}\left(\frac{\xi-x}{
\lambda t^{\frac{\al}{2}}}\right)\right] f(\xi )d\xi  - 
$$
\begin{equation}{\label{d/dx(int x_infty)}}
\hspace{5.5cm}-\displaystyle\lim_{\xi \searrow x}\left[\MM_{\al/2}\left(\frac{x+\xi}{
\lambda t^{\frac{\al}{2}}}\right) + \MM_{\al/2}\left(\frac{\xi-x}{
\lambda t^{\frac{\al}{2}}}\right)\right]f(\xi ) .
\end{equation}
\noindent Note that 

$$\hspace{-6.5cm} \int^{x}_{0}\frac{\partial}{\partial x} \left[\MM_{\al/2}\left(\frac{x+\xi}{
\lambda t^{\frac{\al}{2}}}\right) + \MM_{\al/2}\left(\frac{x-\xi}{
\lambda t^{\frac{\al}{2}}}\right)\right] f(\xi )d\xi=     $$
$$\hspace{2cm} =-\frac{1}{\lambda t^{\al/2}} \int^{x}_{0}\left[\WW\left(-\frac{x+\xi}{\lambda t^{\frac{\al}{2}}},-\frac{\al}{2},1-\al\right) + \WW\left(-\frac{x-\xi}{\lambda t^{\frac{\al}{2}}},-\frac{\al}{2},1-\al\right)\right]f(\xi ) d\xi. $$

\noindent Applying Mean Value Theorem, it yields that

$$\hspace{-3.5cm} \left|\int^{x}_{0}\left[\WW\left(-\frac{x+\xi}{
\lambda t^{\frac{\al}{2}}},-\frac{\al}{2},1-\al\right) + \WW\left(-\frac{x-\xi}{\lambda t^{\frac{\al}{2}}},-\frac{\al}{2},1-\al\right)\right]f(\xi ) d\xi\right|= $$
$$\hspace{2.3cm} = \left|\left[\WW\left(-\frac{x+\theta}{
\lambda t^{\frac{\al}{2}}},-\frac{\al}{2},1-\al\right) + \WW\left(-\frac{x-\theta}{\lambda t^{\frac{\al}{2}}},-\frac{\al}{2},1-\al\right)\right]f(\theta )x\right| \quad \theta \in [0,x]. $$
Then, 

\begin{equation}\label{c_x(0,t)-1}
\displaystyle\lim_{x\searrow 0} \int^{x}_{0}\frac{\partial}{\partial x} \left[\MM_{\al/2}\left(\frac{x+\xi}{
\lambda t^{\frac{\al}{2}}}\right) + \MM_{\al/2}\left(\frac{x-\xi}{
\lambda t^{\frac{\al}{2}}}\right)\right] f(\xi )d\xi=  0.
\end{equation}

\noindent On the other side, 

$$\hspace{-6.3cm} \int^{\infty}_{x}\frac{\partial}{\partial x} \left[\MM_{\al/2}\left(\frac{x+\xi}{
\lambda t^{\frac{\al}{2}}}\right) + \MM_{\al/2}\left(\frac{\xi-x}{
\lambda t^{\frac{\al}{2}}}\right)\right] f(\xi )d\xi=      $$
$$\hspace{-1.8cm}= \int^{\infty}_{x}-\frac{1}{\lambda t^{\al/2}}\left[\WW\left(-\frac{x+\xi}{\lambda t^{\frac{\al}{2}}},-\frac{\al}{2},1-\al\right) - \WW\left(-\frac{\xi-x}{\lambda t^{\frac{\al}{2}}},-\frac{\al}{2},1-\al\right)\right]f(\xi)d\xi =$$

$$\hspace{-1.8cm}= \int^{\infty}_{0}-\frac{1}{\lambda t^{\al/2}}\left[\WW\left(-\frac{x+\xi}{\lambda t^{\frac{\al}{2}}},-\frac{\al}{2},1-\al\right) - \WW\left(-\frac{\xi-x}{\lambda t^{\frac{\al}{2}}},-\frac{\al}{2},1-\al\right)\right]f(\xi )d\xi - $$
$$\hspace{2cm}-
 \int^{x}_{0}-\frac{1}{\lambda t^{\al/2}}\left[\WW\left(-\frac{x+\xi}{\lambda t^{\frac{\al}{2}}},-\frac{\al}{2},1-\al\right) - \WW\left(-\frac{\xi-x}{\lambda t^{\frac{\al}{2}}},-\frac{\al}{2},1-\al\right)\right]f(\xi)d\xi.$$         

\noindent Applying Lebesgue Convergence Theorem to the first integral,
$$\hspace{-1.5cm} \displaystyle\lim_{x\searrow 0} \int^{\infty}_{0}-\frac{1}{\lambda t^{\al/2}}\left[\WW\left(-\frac{x+\xi}{\lambda t^{\frac{\al}{2}}},-\frac{\al}{2},1-\al\right) - \WW\left(-\frac{\xi-x}{\lambda t^{\frac{\al}{2}}},-\frac{\al}{2},1-\al\right)\right]f(\xi )d\xi =$$
$$\hspace{-1.3cm}= \int^{\infty}_{0}\displaystyle\lim_{x\searrow 0}-\frac{1}{\lambda t^{\al/2}}\left[\WW\left(-\frac{x+\xi}{\lambda t^{\frac{\al}{2}}},-\frac{\al}{2},1-\al\right) - \WW\left(-\frac{\xi-x}{\lambda t^{\frac{\al}{2}}},-\frac{\al}{2},1-\al\right)\right]f(\xi)d\xi  =$$
$$\hspace{-1.8cm}=\int^{\infty}_{0}-\frac{1}{\lambda t^{\al/2}}\left[\WW\left(-\frac{\xi}{\lambda t^{\frac{\al}{2}}},-\frac{\al}{2},1-\al\right) - \WW\left(-\frac{\xi}{\lambda t^{\frac{\al}{2}}},-\frac{\al}{2},1-\al\right)\right]f(\xi )d\xi =0. $$
%Luego 
%\begin{equation}\label{(2)} =
  %\displaystyle\lim_{x\searrow 0} \int^{\infty}_{0}-\frac{1}{\lambda t^{\al/2}}\left[W\left(-\frac{x+\xi}{\lambda t^{\frac{\al}{2}}},-\frac{\al}{2},1-\al\right) - W\left(-\frac{\xi-x}{\lambda t^{\frac{\al}{2}}},-\frac{\al}{2},1-\al\right)\right]f(\xi )d\xi =   0  . \end{equation}
\noindent For the second integral we apply Mean Value Theorem as before. Then,

\begin{equation}\label{c_x(0,t)-2}
\displaystyle\lim_{x\searrow 0}\int^{\infty}_{x}\frac{\partial}{\partial x} \left[\MM_{\al/2}\left(\frac{x+\xi}{
\lambda t^{\frac{\al}{2}}}\right) + \MM_{\al/2}\left(\frac{\xi-x}{
\lambda t^{\frac{\al}{2}}}\right)\right] f(\xi )d\xi=0.
\end{equation}

\noindent From (\ref{c_x(0,t)-1}) and (\ref{c_x(0,t)-2}), 
$$\dfrac{\partial c_3}{\partial x}(0,t)=0 .
\quad  \quad  \quad \quad\quad  \quad  \quad $$

\begin{theorem}\label{sol PVIF-IC-nula} Let be $f $ a continuous bounded function in $\bbR^+_0$, then
$$
c(x,t) =\frac{1}{2\lambda t^{\frac{\al}{2}}}\int^{\infty}_{0}\left[\MM_{\al/2}\left(\frac{x+\xi}
{\lambda t^{\frac{\al}{2}}}\right) + \MM_{\al/2}\left(\frac{\left|x-\xi\right|}{
\lambda t^{\frac{\al}{2}}}\right)\right]f(\xi )
 d\xi  \quad \quad \quad \quad \,
$$
is a solution to problem
$$\left\{
\begin{array}{lll}
D^{\al} c(x,t)=\lambda^2\dfrac{\partial^2c }{\partial x^2}(x,t), & &  0<x<\infty, \, 0<t<T, \, 0<\al<1, \\

          \dfrac{\partial c }{\partial x}(0,t)=0,  & &    0<t<T, \\
           c(x,0)=f(x),   & & 0<x<\infty. \\
           %\displaystyle\lim_{x\rightarrow +\infty}c(x,t)=0 & & t>0
\end{array}\right. $$
\end{theorem}

\section{The Initial--Boundary--Value Problem for the Time--Fractional Diffusion Equation in the Quarter Plane with Flux--Boundary Condition.}

\noindent In this last section,  the following problem will be solved:

\begin{equation}{\label{PVIF-IC}}
\left\{\begin{array}{lll}
         _0 D^{\alpha}_t c(x,t)=\lambda^2\dfrac{\partial^2c }{\partial x^2}(x,t) &  & 0<x<\infty, \, 0<t<T, \, 0<\al<1 \\

          c(x,0)=f(x) & \quad  \quad  & 0<x<\infty \\

          \frac{\partial}{\partial x}c(0,t)=g(t) & \quad  \quad  & 0<t<T \\

           %\displaystyle\lim_{x\rightarrow +\infty}c(x,t)=0 &\quad  \quad & t>0
                                       \end{array}\right.  \end{equation}

\noindent As in the previous sections  two auxiliary problems are considered:
\begin{equation}{\label{PVIF-ICA}}
\left\{\begin{array}{lll}
         _0 D^{\alpha}_t c_4(x,t)=\lambda^2\dfrac{\partial^2c_1 }{\partial x^2}(x,t) &  & 0<x<\infty, \, 0<t<T, \, 0<\al<1 \\

          c_4(x,0)=f(x) & \quad  \quad  & 0<x<\infty \\

          \frac{\partial}{\partial x}c_4(0,t)=0 & \quad  \quad  & 0<t<T \\
           %\displaystyle\lim_{x\rightarrow +\infty}c(x,t)=0 &\quad  \quad & t>0
        \end{array}\right.\end{equation}
and

\begin{equation}{\label{PVIF-ICB}}
\left\{\begin{array}{lll}
         _0 D^{\alpha}_t c_5(x,t)=\lambda^2\dfrac{\partial^2c_2 }{\partial x^2}(x,t) &  & 0<x<\infty, \, 0<t<T, \, 0<\al<1 \\

          c_5(x,0)=0 & \quad  \quad  & 0<x<\infty \\

          \frac{\partial}{\partial x}c_5(0,t)=g(t) & \quad  \quad  & 0<t<T \\

           %\displaystyle\lim_{x\rightarrow +\infty}c(x,t)=0 &\quad  \quad & t>0
                                       \end{array}\right.\end{equation}

\noindent From Theorem \ref{sol PVIF-IC-nula},
$$
c_4(x,t) =\frac{1}{2\lambda t^{\frac{\al}{2}}}\int^{\infty}_{0}\left[M_{\frac{\al}{2}}\left(\frac{x+\xi}
{\lambda t^{\frac{\al}{2}}}\right) + M_{\frac{\al}{2}}\left(\frac{\left|x-\xi\right|}{
\lambda t^{\frac{\al}{2}}}\right)\right]f(\xi )
 d\xi  \quad \quad \quad \quad \,
$$

\noindent is a solution to (\ref{PVIF-ICA}).

\noindent In view of the results obtained  in Section 3,
\begin{equation}\label{sol PVIF-ICB} c_5(x,t)=- \int^{\infty}_{x}\int^{t}_{0}W\left(\frac{-\xi}{\lambda (t-\tau)^{\al/2}},-\frac{\al}{2},1\right) g'(\tau) d\tau d\xi
\end{equation}
is proposed as   a solution to problem (\ref{PVIF-ICB}).

\begin{lemma}\label{int_x c es sol} Let $c(x,t)$ be a solution of the time--fractional diffusion equation $_0D^{\al}_t c(x,t)=\lambda^2\dfrac{\partial^2c}{\partial x^2}(x,t)$ such that:
\begin{equation}\label{hip lema integral 1}
\hspace{-1.5 cm}(i) \quad \text{ For every } (x,t),\text{ the function } F(x,t)=\int^{\infty}_{x}c(\xi,t)d\xi  \text{ is well defined, }
  \end{equation}
 \begin{equation}\label{hip lema integral 2}
\hspace{-1 cm}(ii) \quad \displaystyle\lim_{x \rightarrow  \infty}\dfrac{\partial c}{\partial x}(x,t)=0\, , \hspace{9 cm}
 \end{equation}
\begin{equation}\label{hip lema integral 3}
\hspace{-1 cm}(iii) \quad  \left|\frac{\partial}{\partial \tau}c(\xi,\tau)\right|\leq g(\xi) \in L^1(x,\infty)\, , \hspace{7 cm}
\end{equation}
\begin{equation}\label{hip lema integral 4}
\hspace{-1 cm} (iv)  \quad \frac{\frac{\partial}{\partial \tau}c(\xi,\tau)}{(t-\tau)^\al} \in L^1((x,\infty)\times(0,t))\, . \hspace{7 cm}
 \end{equation}

\noindent Then $\int^{\infty}_{x}c(\xi,t)d\xi$ is a solution to the time fractional diffusion equation.
\end{lemma}

\noindent \textit{Proof.}\\

  \noindent From (\ref{hip lema integral 1}), $F(x,t)=\int^{\infty}_{x}c(\xi,t)d\xi$ is well defined. The next equalities are valid due to (\ref{hip lema integral 2}), (\ref{hip lema integral 3})  and (\ref{hip lema integral 4}):
$$_0 D^{\alpha}_t F(x,t)=\frac{1}{\G(1-\al)}\int^{t}_{0}\frac{\partial/\partial \tau F(x,\tau)}{(t-\tau)^\al}d\tau=\frac{1}{\G(1-\al)}\int^{t}_{0}\frac{1}{(t-\tau)^\al}\left(\frac{\partial}{\partial \tau}\int^{\infty}_{x}c(\xi,\tau)d\xi\right)d\tau  = $$
$$\hspace{1.9cm} =\frac{1}{\G(1-\al)}\int^{t}_{0}\frac{1}{(t-\tau)^\al}\int^{\infty}_{x}\frac{\partial}{\partial \tau}c(\xi,\tau)d\xi d\tau=\int^{\infty}_{x}\frac{1}{\G(1-\al)}\int^{t}_{0}\frac{\partial/\partial \tau c(\xi,\tau)}{(t-\tau)^\al}d\tau =$$
$$\hspace{-1cm}=\int^{\infty}_{x}\, _0 D^{\alpha}_t c(\xi,t)d\xi=\int^{\infty}_{x} \lambda^2\dfrac{\partial^2c}{\partial x^2}(\xi,t)d\xi=
-\lambda^2\left.\dfrac{\partial c}{\partial x}(\xi,t)\right|^{x}_{\infty}= $$
$$\hspace{1.4cm}= -\lambda^2 \dfrac{\partial c}{\partial x}(x,t)+\lambda^2 \displaystyle\lim_{x \rightarrow  \infty}\dfrac{\partial c}{\partial x}(x,t)= \dfrac{\partial^2}{\partial x^2}
\left(\int^{\infty}_{x}c(\xi,t)d\xi\right)=\dfrac{\partial^2}{\partial x^2}F(x,t) . \quad \quad $$

\qed

\noindent It can be proved that (\ref{sol PVIF-ICB}) is under the hypothesis of Lemma \ref{int_x c es sol}.\\

\noindent Respect on the border conditions:\\

\noindent $\bullet$	Observing that 
	$$ \int_0^t \MM_{\al/2}\left(\frac{\xi}{\lambda (t-\tau)^{\al/2}}\right)\frac{\xi}{\lambda (t-\tau)^{\al/2+1}}\frac{\al}{2}g(\tau)d\tau\leq
M  W\left(-\frac{\xi}{\lambda t^{\al/2}},-\frac{\al}{2},1\right),  $$
\noindent Lebesgue Convergence Theorem can be applied and  

$$\hspace{-3cm}c_5(x,0)= \lim_{t\searrow 0} - \int^{\infty}_{x} \int_0^t \MM_{\al/2}\left(\frac{x}{\lambda (t-\tau)^{\al/2}}\right)\frac{x}{\lambda (t-\tau)^{\al/2+1}}\frac{\al}{2}g(\tau)d\tau d\xi.
 $$
$$\hspace{-1.6cm}= \lim_{t\searrow 0} \left| \int^{\infty}_{x} \int_0^t \MM_{\al/2}\left(\frac{\xi}{\lambda (t-\tau)^{\al/2}}\right)\frac{\xi}{\lambda (t-\tau)^{\al/2+1}}\frac{\al}{2}g(\tau)d\tau d\xi \right|=  $$
 	
	$$\hspace{-1.6cm} = \left| \int^{\infty}_{x}\lim_{t\searrow 0}\int_0^t \MM_{\al/2}\left(\frac{\xi}{\lambda (t-\tau)^{\al/2}}\right)\frac{\xi}{\lambda (t-\tau)^{\al/2+1}}\frac{\al}{2}g(\tau)d\tau d\xi  \right|\leq $$
	
	$$ \hspace{-5cm}\leq  \int^{\infty}_{x}\left|\lim_{t\searrow 0}M  \WW\left(-\frac{\xi}{\lambda t^{\al/2}},-\frac{\al}{2},1\right)\right|d\xi =0  $$

\noindent	$\bullet$ From  (\ref{cond de borde PVIC x-0}),
	$$\hspace{-1.5cm} \frac{\partial}{\partial x}c_5(0,t)=\lim_{x\searrow 0}\frac{\partial}{\partial x}\left(- \int^{\infty}_{x} \int_0^t \MM_{\al/2}\left(\frac{\xi}{\lambda (t-\tau)^{\al/2}}\right)\frac{\xi}{\lambda (t-\tau)^{\al/2+1}}\frac{\al}{2}g(\tau)d\tau d\xi\right)= $$
	$$\hspace{-1.3cm} =\lim_{x\searrow 0}\int_0^t \MM_{\al/2}\left(\frac{x}{\lambda (t-\tau)^{\al/2}}\right)\frac{x}{\lambda (t-\tau)^{\al/2+1}}\frac{\al}{2}g(\tau)d\tau= g(t).$$

\begin{theorem}\label{sol PVIF-IC-FLUJO}
Let be $f $ a continuous bounded function in $\bbR^+_0$ and $g$ a continuous function in  $[0,T)$. Then 

$$\hspace{-1.4cm}c(x,t)=\frac{1}{2\lambda t^{\frac{\al}{2}}}\int^{\infty}_{0}\left[\MM_{\al/2}\left(\frac{x+\xi}
{\lambda t^{\frac{\al}{2}}}\right) + \MM_{\al/2}\left(\frac{\left|x-\xi\right|}{
\lambda t^{\frac{\al}{2}}}\right)\right]f(\xi )-\quad \quad \quad \quad \quad \quad \quad \quad$$
\begin{equation}\label{ sol PVIF-IC}
\quad\quad \quad\quad \quad \quad \quad\quad \quad \quad  -  \int^{\infty}_{x} \int_0^t \MM_{\al/2}\left(\frac{x}{\lambda (t-\tau)^{\al/2}}\right)\frac{x}{\lambda (t-\tau)^{\al/2+1}}\frac{\al}{2}g(\tau)d\tau d\xi
\end{equation}

\noindent is a solution to problem
\begin{equation}\label{PVIC-FLUJO}
\left\{\begin{array}{lll}
          D^{\al} c(x,t)=\lambda^2\dfrac{\partial^2c }{\partial x^2}(x,t), &  & 0<x<\infty, \, 0<t<T, \, 0<\al<1, \\

          c(x,0)=f(x), & \quad  \quad  & 0<x<\infty, \\

          \frac{\partial}{\partial x}c(0,t)=g(t), & \quad  \quad  & 0<t<T. \end{array}\right.   \end{equation}

\end{theorem}

\begin{theorem}The limit when $\al\nearrow 1$ of the solution to problem 
\begin{equation}\label{PVIC-IC-FLUJO-fnula}\left\{\begin{array}{lll}
          D^{\al} c_5(x,t)=\lambda^2\dfrac{\partial^2c_2 }{\partial x^2}(x,t), &  & 0<x<\infty, \, 0<t<T, \, 0<\al<1, \\

          c_5(x,0)=0, & \quad  \quad  & 0<x<\infty, \\

          \frac{\partial}{\partial x}c_5(0,t)=g(t), & \quad  \quad  & 0<t<T, \\

                                       \end{array}\right. %\end{equation}
\end{equation} is the classical solution to the analogous  problem when $\al=1$ and we recover the heat equation
\begin{equation}{\label{PVI-FLUJO ec del calor}}
\left\{\begin{array}{lll}
          \frac{\partial}{\partial t}w(x,t)=\lambda^2\dfrac{\partial^2c }{\partial x^2}(x,t), &  & 0<x<\infty, \, 0<t<T, \\

          w(x,0)=0, & \quad  \quad  & 0<x<\infty, \\

          \frac{\partial}{\partial x}w(0,t)=g(t), & \quad  \quad  & 0<t<T. \\

                                       \end{array}\right.\end{equation}
\end{theorem}
\proof
It can be seen in \cite{Cannon} that 
\begin{equation}\label{sol-PVI-FLUJO ec del calor}
w(x,t)= - \int^{t}_{0}\frac{e^{-\frac{x^2}{4(t-\tau)}}}{\sqrt{\pi(t-\tau)}}g(\tau)d\tau
 \end{equation}
is a solution to problem (\ref{PVI-FLUJO ec del calor}).\\
 
\noindent Let be $c^{\al}_5$ a solution to problem (\ref{PVIC-IC-FLUJO-fnula}) given by  Theorem \ref{sol PVIF-IC-FLUJO}, 
$$ c^{\al}_5(x,t)=  -  \int^{\infty}_{x} \int_0^t \MM_{\al/2}\left(\frac{x}{\lambda (t-\tau)^{\al/2}}\right)\frac{x}{\lambda (t-\tau)^{\al/2+1}}\frac{\al}{2}g(\tau)d\tau d\xi .$$

\noindent Applying Lebesgue Convergence Theorem, Lemma \ref{conv M al/2 cuando al tiende a 1} and Fubini's Theorem, 

$$\hspace{-1.4cm}\lim_{\al \nearrow 1 } c^{\al}_5(x,t) =\lim_{\al \nearrow 1 }\left\{ -  \int^{\infty}_{x} \int_0^t \MM_{\al/2}\left(\frac{x}{\lambda (t-\tau)^{\al/2}}\right)\frac{x}{\lambda (t-\tau)^{\al/2+1}}\frac{\al}{2}g(\tau)d\tau d\xi . \right\}=$$
$$\hspace{0cm}=-  \int^{\infty}_{x} \int_0^t \lim_{\al \nearrow 1 } \MM_{\al/2}\left(\frac{x}{\lambda (t-\tau)^{\al/2}}\right)\frac{x}{\lambda (t-\tau)^{\al/2+1}}\frac{\al}{2}g(\tau)d\tau d\xi
 = $$
$$\hspace{-4cm}=-\int^{\infty}_{x} \int^{t}_{0}\frac{e^{-\frac{x^2}{4(t-\tau)}}}{2 \sqrt{ \pi}}\frac{x}{\lambda (t-\tau)^{3/2}}g(\tau)d\tau =$$
$$\hspace{1.5cm}=-  \int^{t}_{0} \int^{\infty}_{x}\frac{e^{-\frac{x^2}{4(t-\tau)}}}{2 \sqrt{ \pi}}\frac{x}{\lambda (t-\tau)^{3/2}}g(\tau)d\tau = -\int^{t}_{0}\frac{e^{-\frac{x^2}{4(t-\tau)}}}{\sqrt{\pi(t-\tau)}}g(\tau)d\tau= w(x,t). $$

\endproof

\section{Conclusions}

\noindent %Asymptotic behavior of the Wright functions and bounds for the Mainardi and the Wright function $W(-x,\frac{\alpha}{2}, 1)$ in $\mathbb{R}^+$ have been studied. Then, the analysis of the solution of the initial-value problem for the time-fractional diffusion equation was done, obtaining some conditions for the data to having a solution for the problem. Next,
 On the basis of the asymptotic behavior of some Wright functions and the existence of bounds for the Mainardi and the Wright function $W(-x,\frac{\alpha}{2}, 1)$ in $\mathbb{R}^+$ , three different  initial-boundary value problems for the time-fractional diffusion equation in the quarter plane were solved (considering temperature boundary condition, null flux boundary condition and flux boundary condition in the fixed face $x=0$). In each case, certain conditions must be verified for the data to obtain the solution and the convergence of this solution when $\alpha \nearrow 1$ was analyzed, recovering the classical solutions of the respective boundary-value problems corresponding to the heat equation in the quarter plane.

%%%%%%%%%%%%%%%%%%%%%%%%%%%%%%%%%%%%%%%%%%%%%%%%%
\section{Acknowledgements}

\noindent This paper has
been sponsored by the Project ING495 ``Estudio de diversos problemas con
ecuaciones diferenciales fraccionarias'' from Universidad Nacional de
Rosario, Argentina.

\end{document}